\documentclass[journal]{IEEEtran}
\ifCLASSINFOpdf
\else
\fi
\usepackage{graphicx}
\usepackage{epsfig} 
\usepackage{epstopdf} 
\usepackage{times} 
\usepackage{amsmath} 
\usepackage{amssymb}  
\usepackage{amsfonts}
\usepackage{cite}
\usepackage{multirow, makecell}
\usepackage{subfigure}
\usepackage[utf8]{inputenc}
\DeclareUnicodeCharacter{2212}{-}
\usepackage{color,soul}
\usepackage{tikz}
\usetikzlibrary{positioning}

\usepackage{todonotes}
\presetkeys{todonotes}{inline,backgroundcolor=red!40,caption={}}{}

\usepackage{tablefootnote}
\begin{document}

\title{Stochastic Approximation and Brownian Repulsion based Evolutionary Search}

\author{Rajdeep~Dutta,~T~Venkatesh~Varma,~Saikat~Sarkar,~Mariya~Mamajiwala,~Noor~Awad,~Senthilnath~Jayavelu,~and~Debasish~Roy

\thanks{R. Dutta and S. Jayavelu are with the Institute for Infocomm Research (I$^{2}$R), Agency for Science, Technology and Research (A*STAR), Singapore, 138632.}
\thanks{T.V. Varma is with the Faculty of Mechanical Engineering, Technion–Israel Institute of Technology, Haifa 32000, Israel.}
\thanks{S. Sarkar is with the Department of Civil Engineering, Indian Institute of Technology Delhi, Hauz Khas, New Delhi, 110016, India (corresponding author, e-mail: saikat@iitd.ac.in).}
\thanks{M. Mamajiwala is with the Department of Computer Science, University of Sheffield.}
\thanks{N. Awad is with the Machine Learning Lab, University of Freiburg, Germany.}
\thanks{D. Roy is with the Centre of Excellence in Advanced Mechanics of Materials, Indian Institute of Science, Bangalore, India.}
}

\maketitle

\begin{abstract}
Many global optimization algorithms of the memetic variety rely on some form of stochastic search, and yet they often lack a sound probabilistic basis. Without a recourse to the powerful tools of stochastic calculus, treading the fine balance between exploration and exploitation could be tricky. In this work, we propose an evolutionary algorithm (EA) comprising two types of additive updates. The primary update utilizes stochastic approximation to guide a population of randomly initialized particles towards an optimum. We incorporate derivative-free Robbins-Monro type gains in the first update so as to provide a directional guidance to the candidate solutions. The secondary update leverages stochastic conditioning to apply random perturbations for a controlled yet efficient exploration. Specifically, conceptualized on a change of measures, the perturbation strategy discourages two or more trajectories exploring the same region of the search space. Our optimization algorithm, dubbed as SABRES (Stochastic Approximation and Brownian Repulsion based Evolutionary Search), is applied to CEC-2022 benchmark functions on global optimization. Numerical results are indicative of the potentialities of SABRES in solving a variety of challenging multi-modal, non-separable, and asymmetrical benchmark functions.
\end{abstract}


%
\IEEEpeerreviewmaketitle

\section{Introduction}
Global optimization has widespread applications in science and engineering -- from inverse problems in mechanics, navigation and guidance, automation, optimal control and filtering to integrated circuit design, reconstruction of protein structures and machine learning, and traveling salesman problem (\hspace{-0.1em}\cite{mechopt, reviewOptim, bookGA, SAconv_rate2, ieeeCyber_2024}). At the heart of optimization lies the extremization of an objective function which may be non-convex, non-separable and even non-smooth in real-world applications. This disqualifies gradient-based local optimization strategies to be applicable to several challenging problems -- a reason for the extensive research in designing global optimizers that employ stochastic search via evolutionary algorithms of heuristic or meta-heuristic origin (\hspace{-0.1em}\cite{reviewOptim, algoABC, cuckooSearch1, papDE2}). The usefulness of an optimizer comes from its efficacy in tackling constraints and in scaling across dimensions in the presence of non-convexity. The dimension of an optimization problem depends on the number of decision variables. The constraints on the decision variables arise from physical requirements and regulations, environmental limitations, or modeling approximations. In a bounded search space, it is far more difficult to reach the global optimum of an objective function than just arriving at a local optimum. A traditional gradient-based approach \cite{reviewOptim} needs the derivative information of an objective function to guide the search. Such information quickly directs a candidate solution towards the nearest optimum, although it restricts the search to a local solution only. Even if differentiable, a function's derivative may not always be easy to evaluate or readily available, especially when the objective function represents a one-to-many mapping as encountered in non-convex cases such as many inverse problems \cite{inverse2012} or when it is corrupted with noise \cite{geneticEM}. 

A focus of recent research has therefore been on schemes that avoid derivative calculations (\hspace{-0.1em}\cite{reviewOptim, Saikat1, Saikat2, Bottou1998, Tanabe2020}). In order to overcome the local traps, one may exploit noise, viz. a Monte Carlo method to obtain an empirical distribution involving the decision variables by repeated random sampling of inputs \cite{papSAMC}. Monte Carlo simulations have numerous applications in filtering, estimation, finance, adaptive control and stochastic optimal control \cite{book_Appl}. During the past few decades, evolutionary algorithms (EAs) such as the genetic algorithm (GA) \cite{bookGA}, simulated annealing (SA) \cite{papSA}, differential evolution (DE) \cite{papDE2} and particle swarm optimization (PSO) \cite{papPSO}, have gained prominence thanks to their derivative-free directional guidance features. EAs rely on randomly scattered candidate solutions that mimic nature-inspired mechanisms and use heuristics or meta-heuristics, to incrementally update those candidates over iterations. Such heuristics-based propagators: genetic algorithm \cite{bookGA} based on Darwin's law of the survival of the fittest, artificial bee colony algorithm (\hspace{-0.1em}\cite{algoABC, otherABC}) based on foraging strategies of honey bees, cuckoo search algorithm \cite{cuckooSearch1} based on Levy flight mechanism of cuckoos, have also been used in solving continuous and combinatorial optimization problems. However, towards a principled, efficacious and economical search, these algorithms do not typically exploit the powerful machinery of probability theory or stochastic calculus to inform the directional updates applied to the candidate processes. To be sure, there are a few exceptions; authors in \cite{Saikat1}, for instance, developed a global optimizer, COMBEO, based on a change of measure technique.  

The current work adopts a probabilistic framework as the backbone of the proposed algorithm. We leverage stochastic approximation (SA) 
 and stochastic conditioning through a change of measures to come up with a strategy for global search. The theory of SA is used to give directional guidance to our algorithm. SA is useful in optimizing functions that cannot be evaluated explicitly, but can only be estimated. Estimation of unknown parameters is ubiquitous in many areas, such as (a) system identification problems, (b) adjustable control gains in adaptive control systems, (c) optimal weight matrix in signal processing, and (d) learning parameters in pattern classification. The parameter estimation problem may also be viewed as one of root (zero) finding for an unknown function. Pioneering work in SA was conducted by Robbins and Monro (\hspace{-0.1em}\cite{originalRM, DraftRM}) in the 1950s, and it has since been widely used in myriad applications \cite{bookSA}. Robbins-Monro algorithm has been used across a wide spectrum of problems to sequentially find the zero(s) of a function under noisy or corrupted observations (\hspace{-0.1em}\cite{Kulkarni1, Kulkarni2}). Another SA algorithm by Kiefer-Wolfowitz \cite{stochKW}, known as the finite difference stochastic approximation (FDSA), was used to find the extremum of a function from the available noisy gradient information. A Monte Carlo set-up in combination with the SA has already been used in developing global optimization algorithms for machine learning \cite{papSAMC}. Furthermore, approximate solutions of stochastic optimal control problems \cite{bookSA} can be determined with the help of SA.

An essential key to the success of a global search algorithm lies in its capability of exploring search space efficiently, which fights the risk of getting trapped within the neighborhood of a local optimum \cite{Saikat1}. 
For instance, different mutation techniques are leveraged in genetic algorithm (GA) and matrix adaptation evolution strategy (MA-ES) to induce better exploration \cite{ieeeCyber_2022}. 
The current optimizer draws its exploration power through a stochastic conditioning approach. Stochastic conditioning is a powerful tool with great usefulness in areas such as automatic control, finance, advertising, supply-chain optimization, and dynamic resource allocation \cite{Fleming}. Here, one modifies the dynamics of an evolutionary process under the effects of random influences and control inputs \cite{Note1, condSDE}. It is worth noting that, given the original stochastic process evolving under a known measure, a conditioned stochastic process arises out of an absolutely continuous change of measures  \cite{condBrow, Draft1}. An appropriate change of measures for the underlying sample space has the potential to alter the drift of a stochastic dynamical  system towards certain desirable ends. The present optimizer employs such a change of measures, which is incidentally related to the well known Doob's $h$-transform \cite{Doob1}, to collect information from the unexplored zones of the search space. More specifically, the change of drift may be designed to nudge the stochastic process towards a desired exit state. It is, for example, possible to condition a random walk to always stay positive by applying Doob's $h$-transform \cite{condiBrown}. Similarly, one-dimensional Brownian motions can be conditioned to generate three-dimensional Bessel processes \cite{1dBrown3dBSE}.



According to the no-free lunch theorem \cite{NFL1997}, there is no single algorithm that serves the best in optimizing all types of benchmark functions. Even without the possibility of such an ideal scheme, many real-world applications demand a robust optimizer that can tackle non-smooth, nonlinear, and ill-conditioned problems. The following are the highlights of our contribution in the present work.
\begin{itemize}
\item  SABRES is developed with rigorous reasoning, bypassing dependence on heuristics or meta-heuristics. Arguments based on stochastic approximation and stochastic conditioning, the latter affording a balance of exploitation versus exploration, offer a systematic route to future enhancements in the scheme.  

\item  The directional guidance of SABRES involves Robbins-Monro type gains that are easy to calculate and implement. This makes for computationally inexpensive and easy implementation.


\item An enhanced exploration capability of the proposed algorithm is brought forth by a change of measures, which adds a repulsive drift to the standard Brownian process.

\item The robustness and promise of the algorithm is demonstrated with applications across various multi-modal, non-separable, and asymmetrical benchmark functions.
\end{itemize}

The remainder of the paper is organized as follows. Section \ref{sec:model} introduces the probabilistic tools used in the current work, and Section \ref{sec:method} describes our proposed methodology in combination with the designed operations. In Section \ref{sec:results}, we present numerical results to validate the effectiveness of the proposed optimizer. Finally, Section \ref{sec:conclusion} concludes the work.


\section{The Problem and Our Approach} 
We begin by stating the optimization problem and introducing a stochastic framework to solve it.

\subsection{The Problem} Given a multimodal (non-convex) and possibly non-smooth objective function $f(\mathbf{x}): \Bbb R^D\mapsto \Bbb R$ with $\mathbf{x}\in \Bbb R^{D}$ denoting the $D$-dimensional vector of design variables, the problem is to find the global optimum $\mathbf{x}^*$ such that $f(\mathbf{x})\geq f(\mathbf{x}^*) ~\forall \mathbf{x}$.

\subsection{Stochastic Framework}\label{sec:model}

Consider a complete probability space $(\Omega, \emph{F}, \mathbb{P})$, where the sample space $\Omega$ denotes the set of all elementary events $\omega$ on which the random variables of interest (i.e. the design variables) $\mathbf{x}(\omega):\Omega \rightarrow \Bbb {R}^D$ act. For a less cluttered presentation, we do not notationally distinguish between the mapping $\mathbf{x}(\omega)$ and the values $\mathbf{x}$ taken by $\mathbf{x}(\omega)$ in $\Bbb{R}^D$. The $\sigma$-algebra $\emph{F}$ is constructed with the union, intersection and complementation of all open sets $S\subset \Omega$ and 
 $\mathbb{P}$ is the associated probability measure \cite{DRoyBook}. 
We introduce $\tau$ as a positive monotonically increasing (time-like) variable on $\Bbb R$, useful for tracking the iterations of the dynamic  system we are going to construct. Thus, the design variables are evolved along the iteration axis. 


In this work, we wish to draw upon some aspects of the rich repertoire of diffusive stochastic processes and calculus to develop an appropriate dynamical  system to propagate the design variable, which is considered a vector-valued stochastic process $\mathbf{x}(\tau, \omega)$ or $\mathbf{x}_{\tau}(\omega)$. Thus, for each $\tau$, $\mathbf{x}_{\tau}(\omega)$ is a $\tau$-parametrized, vector-valued random variable. To achieve our optimization goal, we allow that the optimal solution exist within the above probability space. The evolution of the design variable is mimicked by $n$ independently evolving stochastic processes (also referred to as trajectories) along the iteration axis $\tau$, which we denote as $\{ \mathbf{x}_{\tau}^{(1)}(\omega), \mathbf{x}_{\tau}^{(2)}(\omega),...,\mathbf{x}_{\tau}^{(n)}(\omega)  \}$. We emphasize that each trajectory $\mathbf{x}_{\tau}^{(i)}(\omega)$ is a stochastic process in its own right, and not the $i$-th realization of a stochastic process. For further notational simplicity, we may henceforth drop the argument $\omega$ whilst working with these processes. In this stochastic framework, each process $\mathbf{x}_{\tau}^{(i)}~\forall~i \in [1,...,n]$ behaves like a constant-mean martingale \cite{saikat_PLA} once an optimal solution is attained (the mean being the optimal solution). When the stochastic process along $\tau$ converges to a \textit{martingale}, then any future mean of the associated process $\mathbf{x}_{\tau}^{(i)}$ conditioned on the present (on an increasing family of sub-$\sigma$ algebras $N_{\mathcal{T}}$ for $\mathcal{T}\leq \tau$) remains iteration invariant, i.e. $\Bbb E(\mathbf{x}_\tau^{(i)}|N_{\mathcal{T}})=\mathbf{x}_{\mathcal{T}}^{(i)}$, where $N_{\tau}$ denotes the filtration generated by the  processes $\mathbf{x}_{\tau}^{(i)}$ for $i\in [1,n]$.

To start with and in the absence of any other information, we demand that $\mathbf{x}_{\tau}^{(i)}$ be a Brownian motion for each $i$. Once it is appropriately projected (i.e. conditioned) on a so-called extremal cost process (see \cite{saikat_PLA}), $\mathbf{x}_{\tau}^{(i)}$ is expected to approach an optimal solution. The extremal cost process may be identified with the best available objective function value at a given $\tau$. 
Based on such a projection, if we select the best solution at every iteration then $\mathbf{x}_{\tau}^{(i)}~ \rightarrow \mathbf{x}^* ~\forall~i$ as $\tau \rightarrow \infty$, effected by rejection sampling. However, such a search process is typically slow and prone to be trapped by a local optimum. To accelerate the search for the global optimum, we further perturb the evolving design variables. We do so by devising a repulsive drift such that any point (or, rather a small neighbourhood surrounding the point) in the design space is never visited by too many of the $n$ independent stochastic processes. The details are furnished below. 

\noindent \textbf{Prediction:} We first define $n$ i.i.d. Brownian processes, i.e. the stochastic dynamics of the $i^{th}$ process is given by
\begin{eqnarray}\label{eqn:Xpredict}
d \mathbf{x}_{\tau}^{(i)} = \boldsymbol{\gamma} d \mathbf{B}_{\tau}^{(i)}~.
\end{eqnarray}
In (\ref{eqn:Xpredict}), $\mathbf{B}_{\tau}^{(i)} = [B_{1, \tau}^{(i)}, B_{2, \tau}^{(i)},...,B_{D, \tau}^{(i)}]^T \in \Bbb R^{D}$ is a zero-mean $D$-dimensional Brownian motion, and $\boldsymbol{\gamma} \in \Bbb R^{D \times D}$ is the noise intensity matrix with the related covariance matrix being $\boldsymbol{\gamma} \boldsymbol{\gamma}^T \in \Bbb R^{D \times D}$. The states of these independently evolving Brownian processes are stacked in a matrix: $\mathcal{X}_\tau  =[\mathbf{x}_{\tau}^{(1)},\mathbf{x}_{\tau}^{(2)},...,\mathbf{x}_{\tau}^{(n)}] \in \Bbb R^{D \times n}$, and the resulting Brownian process matrix is represented by: $\mathcal{B}_\tau = [\mathbf{B}_{\tau}^{(1)},\mathbf{B}_{\tau}^{(2)},...,\mathbf{B}_{\tau}^{(n)}] \in \Bbb R^{D \times n}$.

As part of our optimization strategy, we solve (\ref{eqn:Xpredict}) within a Monte Carlo (MC) setting to obtain $m$ MC realizations (or particles) for each of the $n$ Brownian processes. Accordingly, ${\bf x}_{\tau}^{(i,k)}=[x_{1, \tau}^{(i,k)}, x_{2, \tau}^{(i,k)},...,x_{D, \tau}^{(i,k)}]^T \in \Bbb R^{D \times 1}$ denotes the $k^{\text{th}} ~ (k\in [1,...,m])$ particle of the $i^{th}$ process, ${\bf x}_\tau^{(i)}$, along the iteration axis $\tau$. 
Within a time-discrete setting, let $\tau_0 < \tau_1 < \tau_2 <...<\tau_l<... < \mathcal{T}$. The $i^{\text{th}}$ discrete random walk form of (\ref{eqn:Xpredict}) is then given by
\begin{equation}\label{eqn:Xpred_disc}
    \mathbf{x}_{\tau_{l+1}}^{-(i,k)} = \mathbf{x}_{\tau_l}^{(i,k)} + \boldsymbol{\gamma} \Delta \mathbf{B}_{\tau_l}^{(i,k)}~;~\Delta \mathbf{B}_{\tau_l}^{(i,k)} = \mathbf{B}_{\tau_{l+1}}^{(i,k)}-\mathbf{B}_{\tau_{l}}^{(i,k)}~,
\end{equation}
where $x_{\tau_{l+1}}^{-(i,k)}$ denotes a prediction of $x_{\tau_l}^{(i,k)}$ at iteration $\tau_{l+1}$.

\noindent \textbf{Additive Updates:}  We propose two types of additive updates to drive the predicted processes to the global optimum. A \textit{directional guidance} step is applied throughout the iterations to steer the candidate solutions (i.e. the individual particles of different trajectories) towards an optimum. A \textit{mutation} step is then applied with low probability towards an enhanced exploration of the search space. The next section explicates on these strategies.

\section{Evolution Strategies}\label{sec:method}
In this section, we discuss the two additive update strategies to provide our algorithm with exploration and exploitation functionalities.  
\begin{figure*}[t]
\centering
  \includegraphics[height=7cm, width=11.5cm]{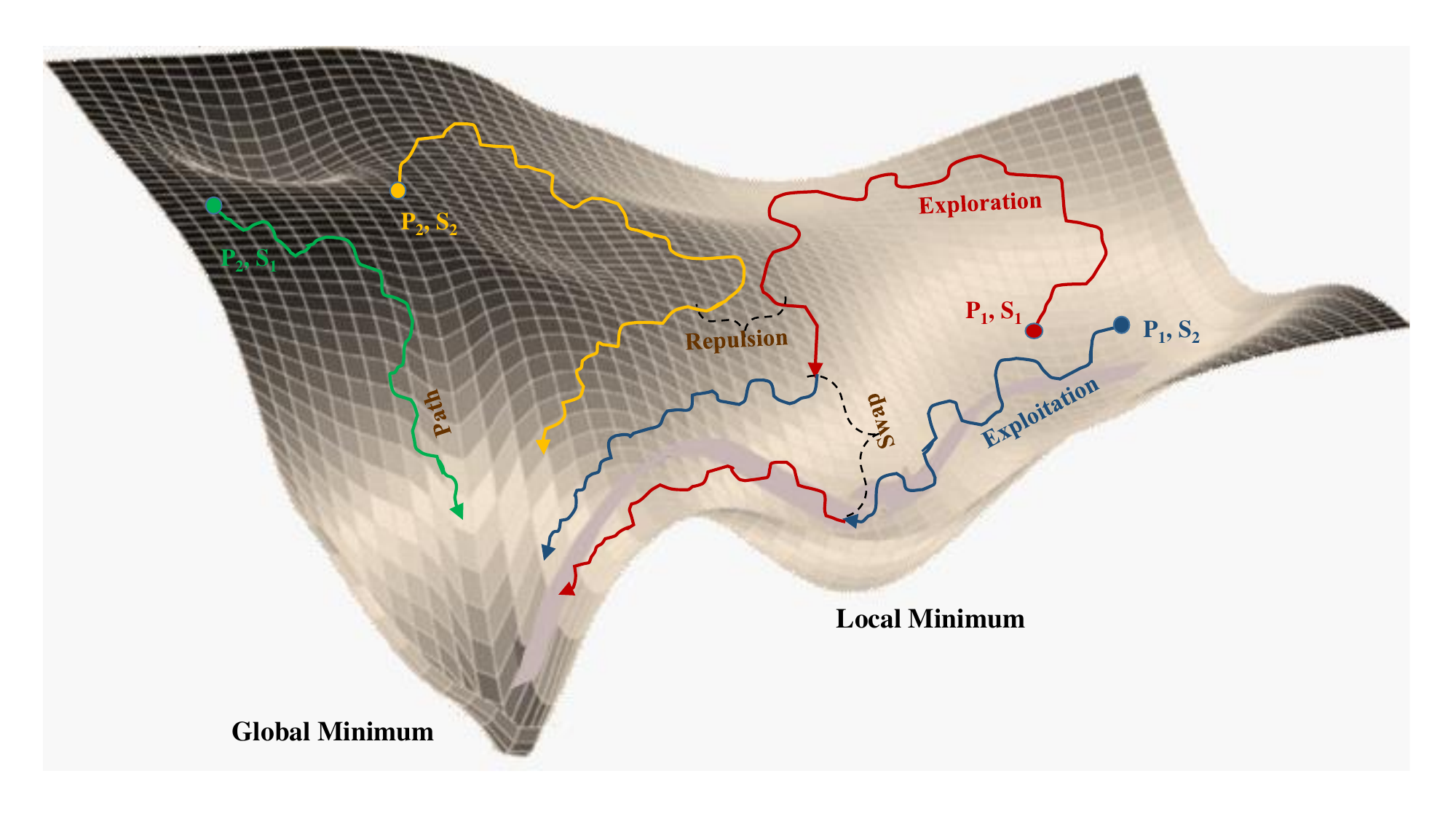}
  \caption{A portrait of the proposed optimization approach: Four sample trajectories are shown on a non-convex function landscape, while considering two particles $P_1,~P_2$ each comprised of two MC samples $S_1,~S_2$. This figure highlights the key concepts like exploration and exploitation involved in the optimization process; exploration includes scrambling (swap) between MC samples and mutual repulsion between different particles.}
  \label{fig:landscape_Optmz}
\end{figure*}

\subsection{Exploration by Stochastic Conditioning}
A well-designed exploration step is imperative to avoid getting stuck at a local optimum. As with most evolutionary algorithms, we propose a mutation step towards an efficient exploration of the search space. However, unlike many other such algorithms, the mutation step we adopt is grounded in stochastic calculus, especially in the Girsanov change of measures. This approach is discussed below. In doing so, we also reflect on how our approach is closely related to Doob's $h$-transform (in its original form) and the non-colliding Brownian motion.

In stochastic calculus, the Girsanov change of measures is used to modify the drift of a stochastic differential equation. However, it is also possible to determine the laws of functionals of diffusion processes based on change of measures \cite{revuzyor}. For instance, given a stochastic process $\mathbf{x}_\tau$ with law $P_x$ (where $\mathbf{x}$ is the initial condition) and a real-valued functional $q(\mathbf{x}_\tau)$ (which could be specified as a constraint or a cost), it is possible to determine a new $q$-mediated process with law $P^q_x$. What is more germane to the present context is that $P^q_x$ could be the law of a process obtained after modifying the original one with the aim of extremizing an appropriate function $q$ or satisfying a constraint involving $q$. In the following, we first describe this concept for a single diffusion process and later demonstrate how we apply this to construct a co-evolving system of multiple (vector valued) diffusions for implementing the mutation step of our algorithm. 

Consider a $D-$dimensional Brownian motion $d\mathbf{x}_\tau = \boldsymbol{\gamma} d\mathbf{W}_\tau$ with generator $\mathcal{L}\phi = \frac{1}{2} c_{ij} \frac{\partial^2 \phi    }{\partial x_i \partial x_j};~~ i,j \in [1,..,D]$, where $\mathbf{c} = \boldsymbol{\gamma} \boldsymbol{\gamma}^T$ and $\boldsymbol{\gamma}$ is the noise intensity matrix. Note that $\mathbf{W}_\tau$ is now a $D$-dimensional vector Brownian motion with standard scalar components (i.e. each scalar component is of mean zero and variance $\tau$). The generator $\mathcal{L}^q$ of the $P^q_x$-process, for a suitable function $q$, is then given in terms of the Carr$\acute{e}$ du Champ operator $\Gamma$ as follows 
\cite{revuzyor} :

\begin{align}
    \mathcal{L}^q \phi &= \mathcal{L} \phi + \Gamma(q, 
\phi),  \\ \nonumber 
\mathrm{where} \;\; \Gamma(q,\phi) &= \mathcal{L}(q\phi) - q \mathcal{L}(\phi) - \phi \mathcal{L}(q),   \\
\mathrm{and}~ \mathcal{L}(q\phi) & = \frac{1}{2} c_{ij} \frac{\partial^2 q \phi    }{\partial x_i \partial x_j} \nonumber \\
&= c_{ij}\frac{\partial q}{\partial x_i} \frac{\partial \phi}{\partial x_j} +   q \mathcal{L}(\phi) + \phi \mathcal{L}(q), \nonumber \\
\mathrm{hence} \;\;
\mathcal{L}^q \phi &= \frac{1}{2} c_{ij} \frac{\partial^2 \phi}{\partial x_i \partial x_j} + c_{ij}\frac{\partial q}{\partial x_i} \frac{\partial \phi}{\partial x_j}.
\end{align}
The corresponding dynamics for the $P^q_x$ process thus becomes 

\begin{equation}
    d\mathbf{x}^{q}_{\tau} = \mathbf{c} \nabla q(\mathbf{x}^{q}_\tau) dt + \boldsymbol{\gamma} d{\mathbf{W}}_\tau .
\end{equation}
Clearly, the modified process is a diffusion that extremizes $q$. 

If we define $r=\exp(q)$, then the modified generator may also be expressed as $\mathcal{L}^r \phi = \mathcal{L} \phi + r^{-1}\Gamma(r,\phi)$. It can be shown that the associated Radon-Nikodym derivative, which is the ratio $dP_x^q/dP_x$, is given by
\begin{equation}\label{eqn:change_of_meas}
    \frac{dP_x^q}{dP_x}=\frac{r(\mathbf{x}_\tau)}{r(\mathbf{x_0})}\exp\left({\int_0^\tau}\frac{\mathcal{L}r(\mathbf{x}_s)}{r(\mathbf{x}_s)}ds\right).
\end{equation}
If we think of $r$ as an unnormalized measure on the sample space $\Omega$, then the change of measure formula in (\ref{eqn:change_of_meas}) guides the diffusion according to $r$. Also note that, when the positive function (the unnormalized measure) $r$ is harmonic, i.e. $\mathcal{L}r=0$, the change of probability measures in  (\ref{eqn:change_of_meas}) pertains exactly to the classical version of Doob's $h$-transform that we exploit in the present work.  

We now direct our attention to an exploitation of this approach for the mutation step. Consider, as before, a system of $n$ independent $D$-dimensional Brownian motions evolving as $d\mathbf{x}^i_\tau = \boldsymbol{\gamma} d\mathbf{B}^i_\tau \;\text{for}\; i \in [1,...,n]$. 
Our objective here is to modify the evolution of these $n$ Brownian trajectories (or processes) so that they together efficiently explore the search space. One possible way to achieve this is to ensure that no two of the $n$ processes ever come too close, thus avoiding the possibility of exploring the same region simultaneously. In organising such a repulsion among the trajectories, we choose the function  $q$ of all $n$ trajectories (each being $D-$dimensional) as follows.
\begin{equation}
q(\mathbf{x}^1,...,\mathbf{x}^n) = \sum_{d=1}^D \sum_{i,j;i<j} \log(x_d^i - x_d^j);~ i,j \in [1,...,n]. 
\end{equation}
where the trajectories are so ordered that the argument of the logarithm always remains positive. In order to appreciate this point better, suppose that they are put in a descending order. In that case, each term $(x_d^i - x_d^j)$ trivially becomes positive since $i<j$. One may then verify that the corresponding function $r=\exp(q)$ is a product of the typical factors $(x_d^i - x_d^j)$ for different $i$ and $j$ and that $r$ is indeed harmonic. Observe that, interpreted as an unnormalized measure, $r$ tends to be zero whenever $x^i_d-x^j_d \rightarrow 0$ for any component $d$ of any two distinct trajectories $\mathbf{x}^i$ and $\mathbf{x}^j$. This is why, a new drift in the form of a repulsive force between two closely separated trajectories arises in the $r$-mediated dynamics. Specifically, the $d^{th}$ element of the additional drift for the $i^{th}$ trajectory is given by
$\frac{\partial q(\mathbf{x}^1,...,\mathbf{x}^n)}{\partial x^i_d} =  \sum_{j, j\neq i} \frac{1}{(x_d^i - x_d^j)}$.
The modified system of $n$ co-evolving $D$-dimensional, drift-modified Brownian trajectories thus obtained is 
\begin{align}\label{eqn:modifiedBM}
dx^i_{d,\tau} = \sum_{j, j\neq i} \frac{d \tau}{\left(x^i_{d,\tau} - x^j_{d,\tau} \right)} + \gamma_{d,d} dB^i_{d,\tau}.
\end{align}


\noindent  \textbf{Exploratory Update:} During implementation, an additive mutation is organized via the following form of equation (\ref{eqn:modifiedBM}), which is applied to the candidate solutions for all components $d\in [1,...,D]$. 
\begin{equation}\label{eqn:update2}
     x_{d, \tau+1}^{-(i,i_\alpha)} = x_{d, \tau}^{(i,i_\alpha)} + \sum_{j, j\neq i, j_\alpha=\tilde{s}(j)} \frac{1}{x_{d, \tau}^{(i,i_\alpha)} - x_{d, \tau}^{(j,j_\alpha)}}+ \gamma_{d,d} \Delta B_{d,\tau}^{(i,i_\alpha)}~, 
\end{equation}
where $x_{d, \tau}^{(i,i_\alpha)}$ denotes the $d^{\text{th}}$ component of the $i_\alpha^{\text{th}}$ MC realization of the $i^{\text{th}}$ particle at iteration $\tau$ and $x_{d, \tau+1}^{-(i,i_\alpha)}$ denotes an exploratory update for $x_{d, \tau}^{(i,i_\alpha)}$ at iteration $\tau+1$ (considering $\Delta \tau=1$). $\Delta B_{d,\tau}^{i,i_\alpha}$ is the associated incremental noise.
Moreover, $i_\alpha =\tilde{s}(i)$, $j_\alpha = \tilde{s}(j)$, where the set $\tilde{s}$ consists of $n$ randomly picked MC realizations so that each of them represents a distinct Brownian motion. Note that $x_{d, \tau+1}^{-(i,i_\alpha)}=x_{d, \tau}^{(i,i_\alpha)}$ when this mutation step is not applied.

It is of interest to observe that Dyson's Brownian motion may be considered as a special case of the more general approach just discussed. For the specific choice of $q$ in this work, where no interaction effects are included and one considers the evolution of one component at a time, the dynamics given by equation (\ref{eqn:modifiedBM}) is indeed the same as Dyson's Brownian motion \cite{Dyson1Brownian_1962}. However, with our present approach, $q$ may be modified to improve upon or tune the repulsive forces at play among the individual trajectories. In this context, it is also worth noting that the modified generator constructed above remains conservative. In the general setting, it is certainly possible to arrive at a non-conservative generator by adding an appropriate source term and thus improve the mutation step.

\subsection{Exploitation of Stochastic Approximation}
The aim is to minimize a function $f(\mathbf{x})$, where $g(\mathbf{x}) \approx \nabla f(\mathbf{x})$ represents the stochastic approximation of its gradient, treated as an unbiased estimator of the true gradient. Such a gradient estimator is known as the \textit{stochastic gradient} \cite{Bottou1998}, which can be calculated numerically using a finite difference between the measured function values at two successive steps.
Computing the stochastic gradient of a function involves less cost than the regular gradient, as it relies on a subset of randomly picked data points instead of the entire set of data points. 
Suppose that $y_{\tau}$ is an unbiased estimator of $\nabla f(\mathbf{x}_{\tau})$ used to update the design variable as  $\mathbf{x}_{\tau+1}=\mathbf{x}_{\tau} - a_{\tau} y_{\tau}$, whilst satisfying: 
\begin{itemize}
    \item that $f$ has a unique minimizer $\mathbf{x}^*$,
    \item that the convexity criterion, viz.  $\inf\limits_{\|\mathbf{x}-\mathbf{x}^* \|_2^2 > a} \nabla f(\mathbf{x})^T(\mathbf{x}-\mathbf{x}^*) >0$, holds $\forall~a>0$ and
    \item that the second moment of the estimator is bounded by  $\mathbb{E}(\|y_{\tau}\|_2^2) \leq A + B\|\mathbf{x}-\mathbf{x}^* \|_2^2$ for some constants $A,~B>0$. 
\end{itemize}
Then the convergence of $\mathbf{x}_{\tau} \rightarrow \mathbf{x}^*$ is ensured subject to the conditions: $a_{\tau} \rightarrow 0$, $\sum_{\tau} a_{\tau} = \infty$ and $\sum_{\tau} a_{\tau}^2 < \infty$.
An explanatory proof of this proposition is provided in \cite{Bottou1998}, which considers a Lyapunov process $L_{\tau}=\| \mathbf{x}_{\tau} -\mathbf{x}^*\|_2^2$ and  shows that the positive variations $L_{\tau}^+=max(0,L_{\tau+1}-L_{\tau})$ satisfy $\sum_{\tau=1}^\infty L_{\tau}^+ < \infty$, implying a  convergent  sequence $L_{\tau}$.

Stochastic approximation (SA) attempts to find the root $\mathbf{x}^*$ of the function $g(\mathbf{x})$. Let $y(\mathbf{x}_{\tau}, \mathbf{\epsilon}_{\tau})=g(\mathbf{x}_{\tau})+\mathbf{\epsilon}_{\tau}$ denote the observed value of $y(\mathbf{x}_{\tau})$ at an instant $\tau$ with $\mathbf{\epsilon}_{\tau}$ being the measurement noise involved. 
The Robbins-Monro algorithm is well established for finding the root (zero) of a function on the Hilbert space $\mathbb{H}$ \cite{Kulkarni1}, where the estimator $\mathbf{x}_{\tau} \in \mathbb{H}$ evolves recursively as 
\begin{eqnarray}\label{eqn:updSA}
\mathbf{x}_{\tau+1}=\mathbf{x}_{\tau+1}^- - a_{\tau} (g(\mathbf{x}_{\tau+1}^-)+\mathbf{\epsilon}_{\tau}) ~.
\end{eqnarray}
Here $a_{\tau}$ is a sequence of positive constants (scalar gains) tending to zero with iterations. Convergence is assured under suitable restrictions on $g$ and the usual assumptions on the gain sequence, i.e. $a_{\tau} \rightarrow 0$, $\sum_{\tau} a_{\tau} = \infty$, $\sum_{\tau} a_{\tau}^2 < \infty$. For the present problem, (\ref{eqn:updSA}) is applied on every process $\mathbf{x}_{\tau}^{(i)}~\forall~i \in [1,...,n]$  with $\mathbf{\epsilon}_{\tau}=\mathbf{0}$ and the gain sequence is chosen as $a_{\tau}=\frac{G_0}{\tau^p}$ with $G_0>0,~0<p \leq 1$.

The above shows that the stochastic gradient of a convex function drives a process to the optimum of the function $f(\mathbf{x})$. 
While global properties can be assessed via local information by optimizing convex functions, this is not true for non-convex functions \cite{sampling_pnas}. 
Towards tackling non-convex problems, an MC set-up promotes an ensemble of trajectories to explore the search space better, which can drive distinct trajectories into different local convex hulls \cite{sampling_pnas}. Further, by conditioning the candidate solutions on an extremal cost process, implemented herein via rejection sampling, the non-convex optimization problem asymptotically approaches a convex optimization problem over iterations \cite{sampling_pnas, Nonconv_1986}, where the theory of stochastic approximation is applicable.

\noindent  \textbf{Directional Guidance:} Throughout the iterative process, each particle evolves in accordance with the following update rule.
\begin{eqnarray}\label{eqn:update1}
\mathbf{x}_{\tau+1}^{(i,k)}=\mathbf{x}_{\tau+1}^{-(i,k)} - \frac{G_0}{\tau^p} \mathbf{I}_{\tau+1}^{-(i,k)}~,
\end{eqnarray}
where $G_0>0$ is a scaling factor and $G_0/\tau^p$ is a decreasing sequence of gains with $p$ being a power index usually chosen so that $0 < p \leq 1$. In (\ref{eqn:update1}), $\mathbf{I}_{\tau}^{-(i,k)} = \textbf{x}_{\tau}^{-(\sigma_n(i),\sigma_m(k))} - \textbf{x}_{\tau}^{-(i,k)}$ is an innovation function, where $\sigma_n(i)$ and $\sigma_m(k)$ denote random permutations over the sets  $\{1,...,n\}$ and $\{1,...,m\}\setminus\{k\}$ respectively \cite{Saikat1}. (\ref{eqn:update1}) helps the particles to get scattered when away from an optimal point, even as it iteratively drives them to the global optimum.

Note that the scrambling step performed through permutations in arriving at the innovation vector in (\ref{eqn:update1}) prevents stalling of the evolution process by applying random swaps between the gain-weighted directional information for distinct trajectories and particles \cite{Saikat1}. In addition to providing a directional guidance, (\ref{eqn:update1}) gives a random perturbation to the evolving trajectories by means of the scrambling operation. 
The major objective of this update is to provide the directional guidance. Although there is some exploration due to the nature of the innovation function $\mathbf{I}$ in (\ref{eqn:update1}), it becomes ineffective when all the trajectories (i.e. processes) begin to converge towards the same local minimum, which is circumvented by our proposed mutation step via \eqref{eqn:update2}. 

\begin{figure}[ht]
\centering
  \includegraphics[scale=0.31]{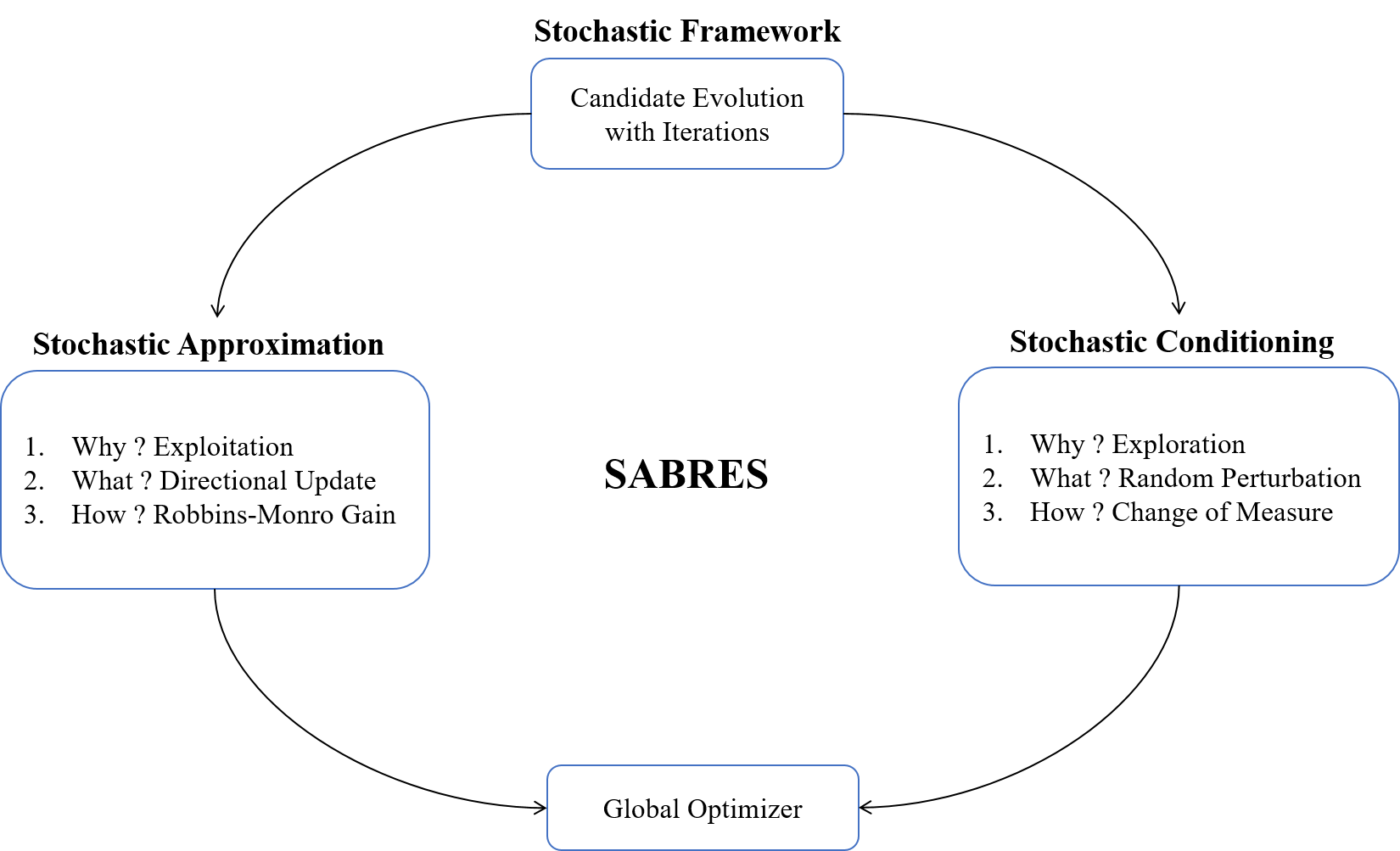}
  \caption{Algorithm design and building blocks: the theoretical basis behind the exploitation and exploration capabilities of SABRES.}
  \label{fig:design}
\end{figure}

\subsection{The Algorithm}
We have so far discussed two major aspects of the proposed optimizer, viz. exploitation and exploration. The basic ingredients of SABRES are highlighted in Figure \ref{fig:design}. 
In the present stochastic framework, novelties are in the theoretical underpinning of the applied random perturbation technique and the simplicity of gains in the directional update. 
At every iteration, an updated solution from the population is selected only if it has a better fitness value than the previous one, otherwise it is discarded. We now provide a pseudocode of the proposed approach followed by a presentation of the related details.

\vspace*{7pt} \hrule \vspace*{1pt}
\noindent \textbf{SABRES Pseudocode:} \vspace*{2pt} \hrule \hrule
{\footnotesize
\begin{itemize}
\item[1.0] -------\textit{Initialization}------- 
\item[] Define a $D$-dimensional bounded search space
\item[] Initialize $n$ i.i.d. Brownian trajectories with $m$ MC realizations for each Brownian trajectory: $\{\{\mathbf{x}_0^{(i,j)}\}_{j=1}^m\}_{i=1}^n$ 
\item[] Initialize parameters - starting gain: $ 5 < G_0 < 20$, gain power: $0 < p\leq 1$, perturbation duration: $\tau_e= 5$, maximum number of iterations: $\mathcal{T}=2 \times 10^3$ or $10^4$.  
\item[2.0]  \textbf{while} $\tau \leq \mathcal{T} $ 
\item[]  ~~\textbf{for} $i=1:n$ 

\item[2.1]  ~~~ -------\textit{Exploration}-------
\item[]  ~~~ \textbf{if}  \textit{Condition $E_r$} (see below) is satisfied 
\item[]  ~~~~~~ Pick one random MC realization $\mathtt{p}_{\alpha}$ from each bin $\mathtt{p}$ 
\item[]  ~~~~~~ to construct the set $\{\mathtt{p}_{\alpha} = \tilde{s}(\mathtt{p}) \}_{\mathtt{p}=1}^n$ 
\item[]  ~~~~~~ \textbf{for} $d=1:D$
\item[]  ~~~~~~~~~ Update particle  position with rule (\ref{eqn:update2}): 
\item[]  ~~~~~~~ $x_{d, {\tau+1}}^{-(i,i_\alpha)} = x_{d, \tau}^{(i,i_\alpha)} + \sum\limits_{j=1,j\neq i}^n \frac{1}{x_{d, {\tau}}^{(i,i_\alpha)} - x_{d, \tau}^{(j,j_\alpha)}} +  \gamma_{d,d} \Delta B_{d,\tau}^{(i,i_\alpha)}$ 
\item[]  ~~~~~~ \textbf{end for} $d$
\item[]  ~~~ \textbf{else}
$x_{d, {\tau+1}}^{-(i,:)} = x_{d, \tau}^{(i,:)}$
\item[]  ~~~ \textbf{end if}

\item[2.2]  ~~~ -------\textit{Exploitation}-------
\item[]  ~~~ \textbf{for} $k=1:m$
\item[]  ~~~~~~~ Construct Innovation: $\mathbf{I}_{\tau+1}^{-(i,k)} = \textbf{x}_{\tau+1}^{-(\sigma_n(i),\sigma_m(k))} - \textbf{x}_{\tau+1}^{-(i,k)}$
\item[] ~~~~~~~ $\sigma_n(i) \in \{1,...,n\}\setminus\{i\},~\sigma_m(k) \in \{1,...,m\}\setminus \{k\}$
\item[]  ~~~~~~~ Multiply Innovation with a decreasing gain: $\frac{G_0}{\tau^p}$
\item[]  ~~~~~~~ \textbf{if}  \textit{Condition $E_t$:} $G_0/\tau^p < 0.05$ is satisfied 
\item[]  ~~~~~~~~~~~ Increase gain value: $G_0=s_g \times G_0;~s_g>0$
\item[]  ~~~~~~~ \textbf{end if}
\item[]  ~~~~~~~ Update particle position with rule (\ref{eqn:update1}):
\item[]  ~~~~~~~ $\mathbf{x}_{\tau+1}^{(i,k)}=\mathbf{x}_{\tau+1}^{-(i,k)} - \frac{G_0}{\tau^p} \mathbf{I}_{\tau+1}^{-(i,k)}$
\item[]  ~~~ \textbf{end for} $k$

\item[2.3]  ~~~ -------\textit{Rejection Sampling}-------
\item[]  ~~~ \textbf{for} $k=1:m$
\item[]  ~~~~~~~~ \textbf{if} $f(\mathbf{x}_{\tau+1}^{(i,k)}) > f(\mathbf{x}_{\tau}^{(i,k)})$ \\
\item[]  ~~~~~~~~~~~ $\mathbf{x}_{\tau+1}^{(i,k)} = \mathbf{x}_{\tau}^{(i,k)}  $  
\item[]  ~~~~~~~~ \textbf{end if}
\item[]  ~~~ \textbf{end for} $k$ 
\item[]  ~~ \textbf{end for} $i$

\item[]  ~~ Update iteration axis: $\tau \rightarrow \tau+1$
\item[]  \textbf{end while}

\item[3.0] \textit{Return} $\{\{\{\mathbf{x}_{\tau}^{(i,j)}\}_{j=1}^m\}_{i=1}^n \}_{\tau = 0}^\mathcal{T}, f(\{\{\{\mathbf{x}_{\tau}^{(i,j)}\}_{j=1}^m\}_{i=1}^n \}_{\tau = 0}^\mathcal{T}) $ 
\end{itemize}
}  
\vspace*{2pt} \hrule \vspace*{7pt}

The convergence characteristics of stochastic approximation algorithms are well established (\hspace{-0.1em}\cite{bookSA, SAEM}). As already noted, we apply SA to provide directional guidance to the evolving candidates. Recall that in the directional update, a sequence of positive gains decreasing with iterations ensures the convergence of SA (\ref{eqn:updSA}). In this context, the selection of an appropriate gain sequence has been a barrier in many applications. Earlier research showed that a sequence going to zero slower than $O(1/n)$ can potentially render convergence at the optimal rate (\hspace{-0.1em}\cite{SAconv_rate1, SAconv_rate2}). Inspired by this, an iterative average technique was proposed in the past to improve the asymptotic convergence rate, while separating the time scales of the averaged and the original processes \cite{SAconv_rate1}. To attain a faster rate of convergence, in SABRES, a slow decline of the gains is ensured by selecting the tuning parameter, $p$, less than one.

The proposed algorithm involves two conditions: \textit{Condition $E_r$} in the exploration phase and \textit{Condition $E_t$} in the exploitation phase. 
\textit{Condition $E_r$} monitors when to apply random perturbations to the candidate solutions. This condition is triggered when there is a small change in the variance of the candidate solutions during their evolution, although the corresponding innovation error is nonzero. 
For example, in our implementation, it is triggered at $\tau$ when var($\{\{f({\bf x}_{\tau-1}^{(i,j)})\}_{j=1}^m\}_{i=1}^n$) $> \frac{1}{2}$ var($\{ \{f({\bf x}_{\tau-2}^{(i,j)})\}_{j=1}^m\}_{i=1}^n\}$) after a regular interval of $1000$ iterations, with a low probability of occurrence (uniform random number $<0.1$). 
Note that perturbations are applied for a few iterations since a prolonged application might impede directional guidance and affect convergence \cite{bound_DysonB}. 
\textit{Condition $E_t$} mitigates the problem of early diminishing gains during the optimization process and monitors when to change the default sequence of decreasing gains. When such a scenario arises first, the primary update gain is restarted with a higher initial value without interrupting its decreasing trend (slope). At the subsequent occurrence(s), the gain is selected randomly from a uniform distribution with an upper threshold. 

\section{Numerical Results}\label{sec:results}
We conduct numerical experiments as per the rules of the IEEE Congress on Evolutionary Computation (CEC) 2022 competition on Single Objective Bound Constrained Real-Parameter Numerical Optimization \cite{benchmark_cec2022}.   
The CEC 2022 benchmark comprises of 12 test functions with dimensions 10 and 20. 
The computational budget or the maximum number of function evaluations ($MaxFEs$) for 10 and 20-dimensional functions are set to 2E+5 and 1E+6, respectively.   
At the beginning of the evolutionary optimization, the search space (bounded by $\pm 100$) is populated with a set of randomly generated (uniformly distributed) candidate solutions.
Note that to maintain reproducibility, the seeds provided by the competition organizers are used for random number generation. 
Also, each trial/run is terminated when either the objective function error value goes below 1E-8 or the computational budget is fully exhausted. 
The selected parameters in SABRES are: number of bins $n=10$, number of MC realizations in each bin $m=10$, initial gain $G_0=10$, gain power $p=0.62$ for $D=10$ and $p=0.7$ for $D=20$.
\begin{table}
    \centering
    \scriptsize
    \caption{Performance statistics of the CEC 2022 test function error values over 30 runs for dimension 10.}
    \begin{tabular}{c|c|ccccc} \hline
        F & Algo & Min & Max & Median & Mean & Std \\ \hline \hline
1 & SABRES & \bf{1.00E-08}  &  \bf{1.00E-08} & \bf{1.00E-08}  & \bf{1.00E-08}  & 1.42E-09 \\
& EA4eig & \bf{1.00E-08} & \bf{1.00E-08} & \bf{1.00E-08} & \bf{1.00E-08} & 1.34E-09 \\ \hline 

2 & SABRES & 1.7019  & 6.2987  & 2.7259 &  3.0651  & 2.0501 \\
& EA4eig & \bf{1.00E-08} & 3.98658 & \bf{1.00E-08} & 1.46175 & 1.95395 \\ \hline

3 & SABRES  & \bf{1.00E-08} &  \bf{1.00E-08}  & \bf{1.00E-08}  & \bf{1.00E-08}  & 1.51E-09 \\
& EA4eig & \bf{1.00E-08} & \bf{1.00E-08} & \bf{1.00E-08} & \bf{1.00E-08} & 9.98E-10 \\ \hline
 
4 & SABRES  & 11.3886 &  23.0693 &  11.6935  &  17.5070 &  30.1106 \\
& EA4eig & \bf{1.00E-08} & 3.97984 & 9.95E-01 & 1.26028 & 1.04298 \\ \hline
 
5 & SABRES  & \bf{1.00E-08}  & 5.91E-05 &  4.29E-07 &  3.99E-06 &  1.14E-05 \\
& EA4eig & \bf{1.00E-08} & \bf{1.00E-08} & \bf{1.00E-08} & \bf{1.00E-08} & 1.62E-09 \\ \hline 
 
6 & SABRES  & 6.03E-01  &  75.9248 &  1.5700 & 5.2807 &  14.2585 \\
& EA4eig & 5.75E-04 & 1.48E-01 & 5.41E-03 & 1.74E-02 & 3.57E-02 \\ \hline

7 & SABRES  & 8.62E-02  & 5.44E-01  & 2.26E-01 &  2.22E-01 &  1.04E-01 \\
& EA4eig & \bf{1.00E-08} & \bf{1.00E-08} & \bf{1.00E-08} & \bf{1.00E-08} & 1.17E-09 \\ \hline

8 & SABRES  & 9.87E-01  &  4.6374  &  1.8043 &  2.0462  & 8.56E-01 \\
& EA4eig & 3.20E-04 & 2.64E-01 & 4.67E-02 & 7.09E-02 & 6.81E-02 \\ \hline

9 & SABRES  & \bf{1.00E-08}  &  359.1124 &  \bf{1.00E-08}  &  12.0896 & 64.4428 \\
& EA4eig & 185.502 & 185.502 & 185.502 & 185.502 & 5.78E-14 \\ \hline

10 & SABRES  & \bf{1.00E-08}  &  \bf{1.00E-08}  & \bf{1.00E-08} & \bf{1.00E-08}  &  212.5104 \\
& EA4eig & 100.084 & 100.213 & 100.1585 & 100.1565 & 3.60E-02 \\ \hline

11 & SABRES  & \bf{1.00E-08}  &  \bf{1.00E-08}  &  \bf{1.00E-08}  &  \bf{1.00E-08}  & 6.92E-10 \\
& EA4eig & \bf{1.00E-08} & \bf{1.00E-08} & \bf{1.00E-08} & \bf{1.00E-08} & 1.06E-09 \\ \hline

12 & SABRES  &  \bf{1.00E-08}  &  150.0581  &  137.6844  &  11.0599  & 66.0403 \\ 
& EA4eig & 145.295 & 158.55 & 145.662 & 147.378 & 3.90145 \\ \hline \hline
    \end{tabular}
    \label{tab:errors_dim10}
\end{table}

The statistics of the error values achieved with SABRES are reported in Tables \ref{tab:errors_dim10} and \ref{tab:errors_dim20}. Considering both the 10 and 20-dimensional problems, SABRES succeeds in solving 14 out of 24 functions overall, which is a little better than the performance (12 out of 24) of the existing top-ranking algorithm, EA4eig \cite{EA4eig}.
SABRES performs well on functions 1, 3, 5, 7, 8, 10, 11, 12, although its performance is not up to the mark for the rest. 
Especially for hybrid and composition functions except 6, SABRES produces promising results. 
The achieved performance for hybrid functions 7 and 8 with 20 dimensions, is the best among all the existing state-of-the-art (SOTA) algorithms.  
SABRES is able to solve composition functions 9 and 12 with 10 dimensions and hybrid function 8 with 20 dimensions, which could not be solved by the existing top-ranking algorithms, i.e., EA4eig \cite{EA4eig}, NL-SHADE-LBC \cite{NL-SHADE-LBC}, and NL-SHADE-RSP-MID \cite{NL-SHADE-RSP-MID}.
For $D=10$, SABRES yields the lowest error in function 12, and for $D=20$, the achieved error value is better than that of NL-SHADE-LBC and NL-SHADE-RSP-MID. 
\begin{figure}%
\centering
\subfigure[][]{%
\label{fig:f9_d10}%
\includegraphics[height=1.5in, width=1.72in]{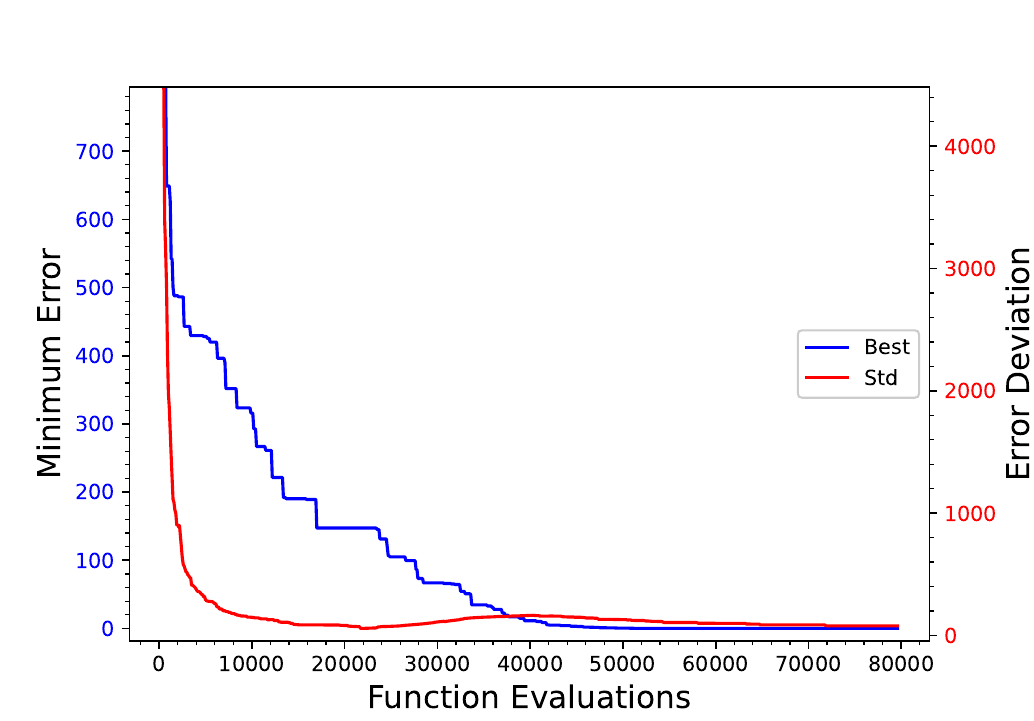}}%
\hspace{2pt}%
\subfigure[][]{%
\label{fig:f12_d10}%
\includegraphics[height=1.5in, width=1.72in]{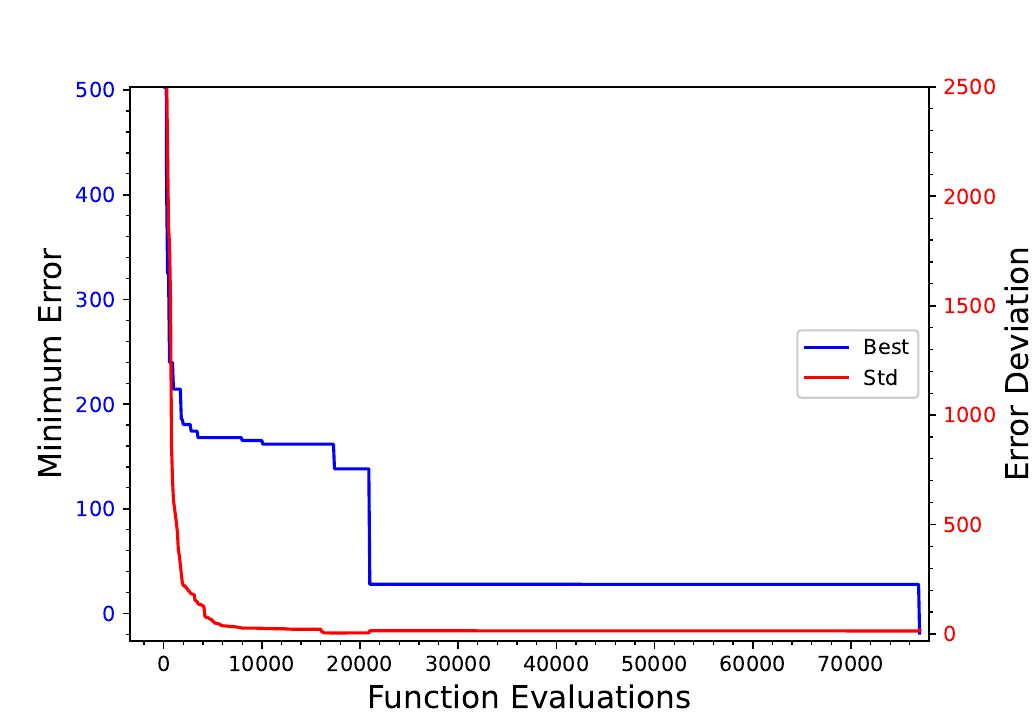}} 
\caption[Evolution of errors for D=10]{Evolution of errors with function evaluations for D=10:
\subref{fig:f9_d10}  the minimum and standard deviation of the error values ($f_9-f_9^*$),
\subref{fig:f12_d10} the minimum and standard deviation of the error values ($f_{12}-f_{12}^*$) of evolving candidate solutions.}%
\label{fig:d10}%
\end{figure}

\begin{figure}%
\centering
\subfigure[][]{%
\label{fig:f7_d20}%
\includegraphics[height=1.58in, width=1.72in]{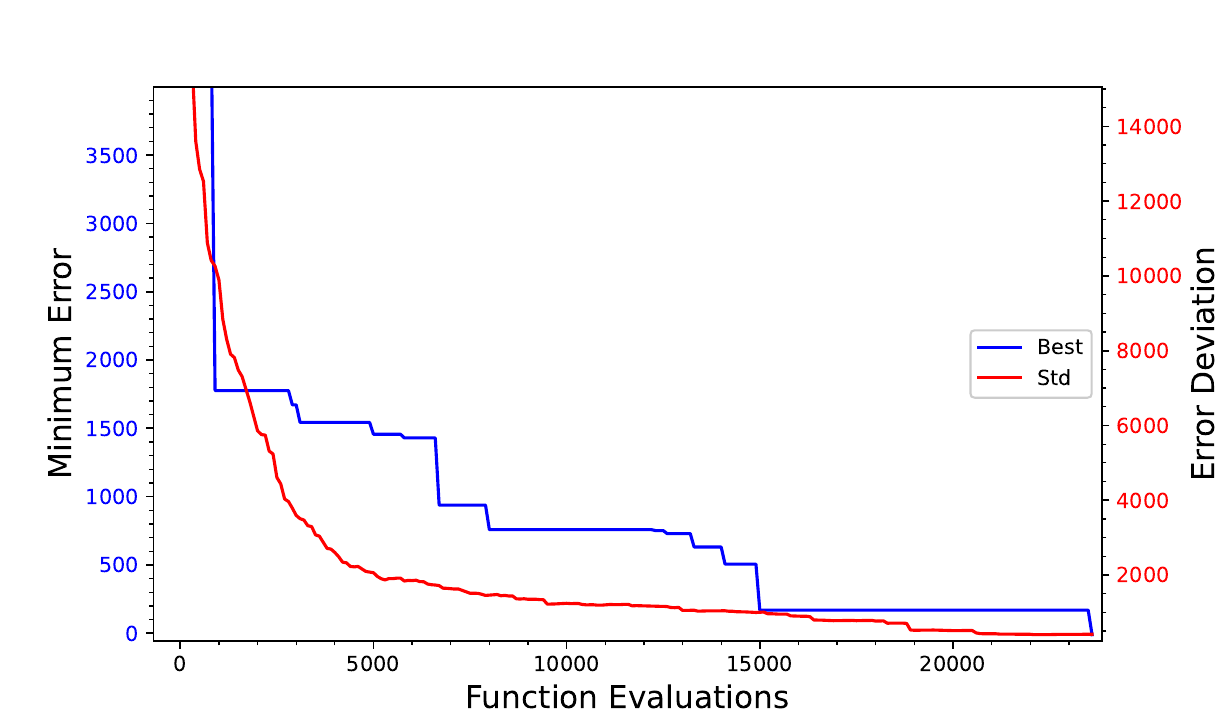}}%
\hspace{2pt}%
\subfigure[][]{%
\label{fig:f8_d20}%
\includegraphics[height=1.58in, width=1.72in]{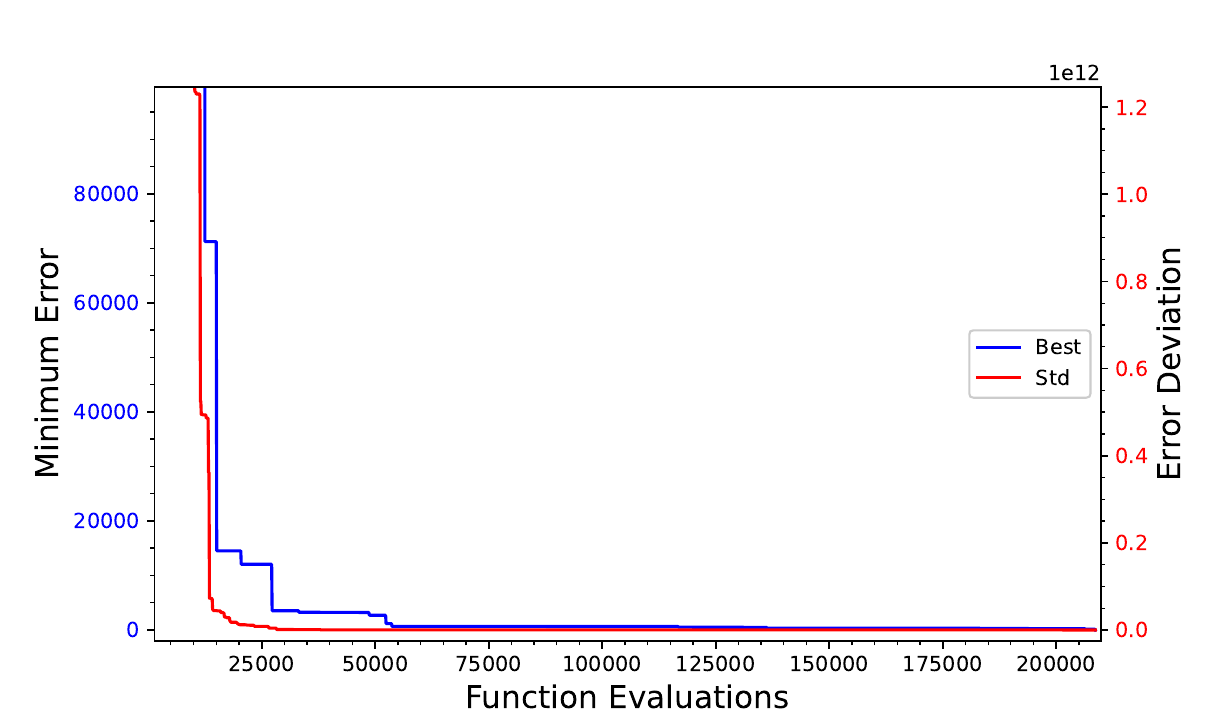}}
\caption[Evolution of errors for D=10]{Evolution of errors with function evaluations for D=20:
\subref{fig:f7_d20} the minimum and standard deviation of the error values ($f_7-f_7^*$),
\subref{fig:f8_d20} the minimum and standard deviation of the error values ($f_{8}-f_{8}^*$) of evolving candidate solutions.}%
\label{fig:d20}%
\end{figure}

Figures \ref{fig:d10} and \ref{fig:d20} showcase a few best-performing case studies, exhibiting how the errors in different objective function values evolve with iterations for $D=10$ and $20$, respectively.
To pay a closer attention to the error evolution, we have zoomed over the axis representing the minimum error.   
In Figure \ref{fig:f9_d10}, there is a noticeable change in the slope of the standard deviation of errors after 3E+4 function evaluations, enabling the candidate solutions to overcome local entrapment (local minimum or saddle point) and arrive at the global optimum.  
This is due to the enhanced exploration capability of SABRES, which facilitates in solving function $9$ that cannot be solved by the SOTA algorithms in 2E+5 evaluations. 
Moreover, SABRES is able to minimize functions $7$ and $8$ efficiently with less function evaluations, as shown in Figures \ref{fig:f7_d20} and \ref{fig:f8_d20}.   
Note that there exist fluctuations in the standard deviation of errors in Figure \ref{fig:f8_d20}, indicating the explorative nature of search. 

\begin{table}
    \centering
    \scriptsize
    \caption{Performance statistics of the CEC 2022 test function error values over 30 runs for dimension 20.}
    \begin{tabular}{c|c|ccccc} \hline
         F & Algo & Min & Max & Median & Mean & Std \\ \hline \hline
1 & SABRES & \bf{1.00E-08}  & 3.84E-04 &  \bf{1.00E-08} &  1.28E-05  & 6.90E-05 \\
& EA4eig & \bf{1.00E-08} & \bf{1.00E-08} & \bf{1.00E-08} & \bf{1.00E-08} & 1.14E-09 \\ \hline

2 & SABRES & 49.0883 &  49.2133 &  49.1263 &  49.1369 &  3.41E-02 \\
& EA4eig & \bf{1.00E-08} & 3.98662 & \bf{1.00E-08} & 1.06310 & 1.79308 \\ \hline

3 & SABRES & \bf{1.00E-08} &  \bf{1.00E-08} &  \bf{1.00E-08} &  \bf{1.00E-08}  & 1.04E-09 \\
& EA4eig & \bf{1.00E-08} & \bf{1.00E-08} & \bf{1.00E-08} & \bf{1.00E-08} & 9.38E-10 \\ \hline

4 & SABRES & 48.3497 &  103.7996 &  73.5883 &  73.9966 &  11.5779 \\
& EA4eig & 3.97984 & 19.8992 & 6.96471 & 8.68931 & 4.08091 \\ \hline

5 & SABRES & \bf{1.00E-08} &  1.10E-04 &  9.52E-07 &  8.75E-06 &  2.06E-05 \\
& EA4eig & \bf{1.00E-08} & \bf{1.00E-08} & \bf{1.00E-08} & \bf{1.00E-08} & 8.89E-10 \\ \hline

6 & SABRES & 9.06E+05 &  2.91E+06 &  1.59E+06 &  1.65E+06 &  4.92E+05 \\
& EA4eig & 2.92E-02 & 4.42E-01 & 1.06E-01 & 1.49E-01 & 1.16E-01 \\ \hline

7 & SABRES & \bf{1.00E-08}  & \bf{1.00E-08}  & \bf{1.00E-08}  &  \bf{1.00E-08}  & 13.0542 \\
& EA4eig & \bf{1.00E-08} & 20 & 2.30209 & 3.50430 & 4.77194 \\ \hline

8 & SABRES & \bf{1.00E-08}  & \bf{1.00E-08} &  \bf{1.00E-08} &  \bf{1.00E-08} &  38.2233 \\
& EA4eig & 2.91E-01 & 21.0441 & 20.2613 & 16.6196 & 7.47147 \\ \hline

9 & SABRES & 335.6402 & 337.9255 &  335.7679  &  336.0313  &  5.40E-01 \\
& EA4eig & 165.344 & 165.344 & 165.344 & 165.344 & 2.89E-14 \\ \hline

10 & SABRES & \bf{1.00E-08} &  \bf{1.00E-08} &  \bf{1.00E-08} &  \bf{1.00E-08} &  340.9627 \\
& EA4eig & 100.201 & 223.252 & 100.256 & 108.259 & 30.4765 \\ \hline

11 & SABRES & \bf{1.00E-08} & 1.14e-01  &  \bf{1.00E-08} &  6.47e-03  & 2.33e-02 \\
& EA4eig & 300 & 400 & 300 & 323.333 & 43.0183 \\ \hline

12 & SABRES & 201.3789 &  219.2431  &  213.4065 &  212.9770 &  4.2520 \\ 
& EA4eig & 188.675 & 200.005 & 200.004 & 199.626 & 2.06839 \\ \hline \hline
    \end{tabular}
    \label{tab:errors_dim20}
\end{table}

\begin{figure}%
\centering
\subfigure[][]{%
\label{fig:example_noE}%
\includegraphics[height=1.54in, width=1.72in]{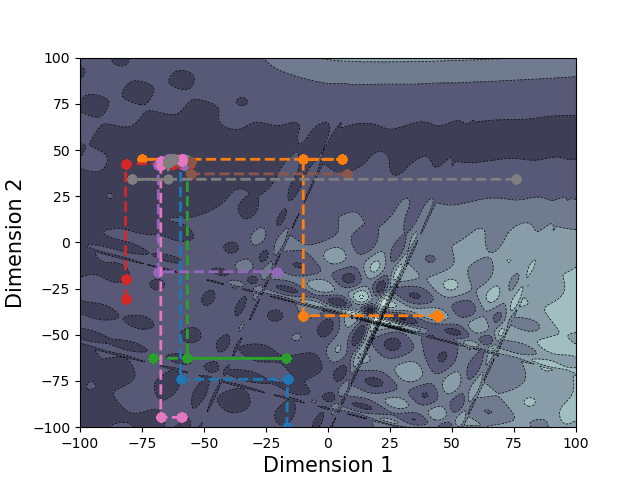}}%
\hspace{2pt}%
\subfigure[][]{%
\label{fig:example_Ee}%
\includegraphics[height=1.54in, width=1.72in]{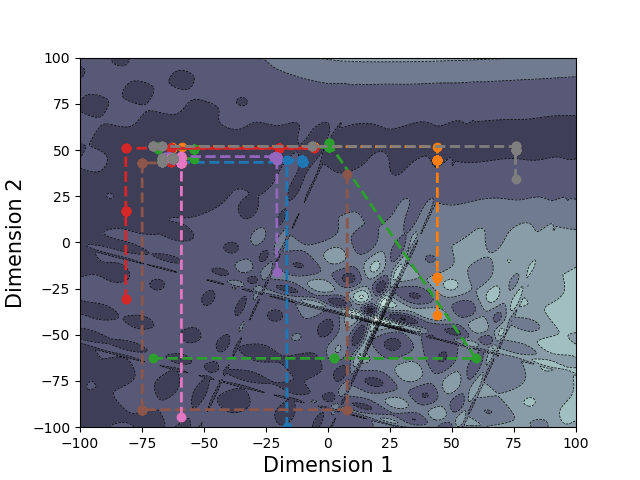}} 
\caption[Paths on Contour for 2D]{2D paths on Contour of F$12$:
\subref{fig:example_noE} evolution of candidate solutions without exploration,
\subref{fig:example_Ee} evolution of candidate solutions with exploration.}%
\label{fig:path2DContour}%
\end{figure}

Further, to gain insight into the exploration of the search space, we consider a simple example of minimizing function $12$ with just $2$ dimensions. 
In this case study, the parameters in SABRES are selected as: number of bins $n=4$, number of MC realizations in each bin $m=2$, initial gain $G_0=1$, gain power $p=0.7$, while retaining the stopping criteria. 
For a consistent initialization of candidate solutions, i.e., fixed initial population owing to a specific seed, we run SABRES with and without the proposed random perturbation technique. 
The resulting trajectories of the candidates on the 2D contour of function 12, are depicted in Figure \ref{fig:path2DContour}.
Here, $8$ candidates evolve iteratively within a search space bounded by ($\pm 100$), and arrive at the global optimum of $[-62.3918, 45.3187]$ finally.   
The 2D paths on the contour reveal that the proposed random perturbation strategy enables candidate solutions to explore the entire search space more efficiently.  
Due to this aspect, SABRES with the exploration capability solves function 12 with the minimum error of 1E-8 in $7235$ iterations ($57888$ evaluations), whereas SABRES without any exploration capability yields the minimum error of $0.0114$ in $25000$ iterations (2E+5 evaluations).  

For performance comparison, we rely on the CEC 2022 ranking scheme that rewards an algorithm not only based on the final objective function error value (accuracy) obtained at the end of a search, but also based on the function evaluations (speed) required to find the optimum.   
It is worth noting that the associated score degenerates to the Mann-Whitney U-statistic in case only two algorithms participate in the competition.    
The best trial in solving each function is identified when the minimum error value of $\leq$ 1E-8 is obtained with the least function evaluations, which is given a score: no. of trials $\times$ no. of algorithms.   
The scores achieved by three algorithms, NL-SHADE-LBC, SABRES, and NL-SHADE-RSP-MID, are presented in Table \ref{tab:rank}.
Clearly, our proposal performs very well on 20D functions and decently on 10D functions, which makes it the $3^{rd}$ (rank 3) algorithm as per the CEC 2022 competition criteria. 
According to Table \ref{tab:rank}, the performance of SABRES on the composition functions for $D=10, 20$ and the hybrid functions for $D=20$, is significantly better than the others.

\begin{table}
    \centering
    \caption{Scores achieved by the top-tier evolutionary algorithms over 30 runs, as per the CEC 2022 ranking criteria. 
    \tablefootnote{The result data available in: https://github.com/P-N-Suganthan/2022-SO-BO, indicate no stopping criteria for the top algorithm, EA4eig, and every trial is run up to the maximum function evaluations. Hence, we have excluded EA4eig from the score table.}}
    \begin{tabular}{c|c|c|c|c|c|c} \cline{1-7}
         Function &  \multicolumn{2}{l|}{NL-SHADE-LBC}  & \multicolumn{2}{c|}{SABRES} & \multicolumn{2}{l}{NL-SHADE-RSP}  \\ \cline{1-7}
         Dimension & 10 & 20 & 10 & 20 & 10 & 20 \\ \hline \hline
         1 & 1800 & 990 & 900 & 1650 & 0 & 60 \\
         2 & 1752 &  959 & 18 &  0 & 930 & 1741 \\
         3 & 0 &  0 & 1800 &  1800 & 900 & 900 \\
         4 & 1798 &  1800 & 87 &  3 & 815 & 897 \\
         5 & 1800 &  1770 & 716 &  930 & 184 & 0 \\
         6 & 1395 &  1760 & 14 &  0 & 1291 & 940 \\
         7 & 900 &  831 & 0 &  1800 & 1800 & 69 \\
         8 & 1485 &  339 & 0 &  1800 & 1215 & 561 \\
         9 & 930 &  1800 & 1740 &  0 & 30 & 900 \\
         10 & 30 &  0 & 1800 &  1800 & 870 & 900 \\
         11 & 1680 &  435 & 150 &  1680 & 870 & 585 \\
         12 & 469 &  694 & 1800 &  1800 & 431 & 206 \\ \hline \hline
         $\sum$ & 14039 &  11378 & 9025 & 13263 & 9336 & 7759 \\ \hline 
    \end{tabular}
    \label{tab:rank}
\end{table}


The algorithmic complexity of SABRES is shown in Table \ref{tab:AlgoComplexGOSAAC}. The associated parameters are as follows: (a) $T_0$ is the running time of a simple computer program comprised of various basic operations and it depends just on the processor in use; (b) $T_1$ is the time consumed in evaluating function $1$ repeatedly (2E+5 times), which accounts for the processor and function evaluation; and (c) $\hat{T}_2$ is the estimated time taken (averaged over $5$ runs) by SABRES to minimize function $1$ for $MaxFEs=$ 2E+5, which accounts for the processor, function evaluation, and algorithmic computation. 
Note that the execution time ($\hat{T}_2$) does not scale up significantly as the underlying problem's dimension increases from 10 to 20. 


\begin{table}[ht]
\centering
\caption{Algorithmic complexity of SABRES}
\label{tab:AlgoComplexGOSAAC}
\begin{tabular}{|c|c|c|c|c|} 
\hline
D & $T_0$ sec & $T_1$ sec & $\hat{T}_2$ sec & $(\frac{\hat{T}_2-T_1}{T_0})$ \\ 
\hline \hline
10 & 0.4998 & 1.2245 & 7.9870 & 13.5304 \\ 
20 & 0.4998 & 1.2416 & 12.5305  & 22.5868 \\
\hline
\end{tabular}
\end{table}

\section{Conclusion}\label{sec:conclusion}
In non-convex, non-smooth optimization, no single algorithm fits all. Nevertheless, the point that we wish to emphasize is that the proposed idea, including the related algorithm - SABRES being the acronym - has a rational grounding. It thus offers a principled means of tuning exploration against exploitation and also a scientifically tenable framework for future enhancements. There are a few desirable features in the present proposal. For instance, gain calculations involved in the primary update are computationally inexpensive. The exploratory aspect, organized through a change of measures and ensuring economy of exploration through repulsion of nearby trajectories, is also novel with a possibility for future enhancements. Accordingly, SABRES proves to be efficient and fast in solving the CEC 2022 benchmark functions. 
In terms of accuracy too, its overall performance is noteworthy - substantively outperforming some of the top-ranking algorithms available in the literature. In future, we plan to improve upon the exploratory aspect of the algorithm by simultaneously requiring that the nearby trajectories repel whilst the far-off ones attract each other. We also intend to extend and exploit SABRES for multi-objective saddle-point problems involving constraints and appearing in myriad practical applications.  

\section*{Acknowledgement}
This research work has been supported by the Department of Civil Engineering and Centre of Excellence in Advanced Mechanics of Materials, Indian Institute of Science (IISc), Bangalore, and the Accelerated Materials Development for Manufacturing Program at A*STAR via the AME Programmatic Fund by the Agency for Science, Technology and Research under Grant No.  A1898b0043.

\end{document}